\documentclass{article}%
\usepackage{amsfonts}
\usepackage{amsmath}
\usepackage{amssymb}
\usepackage{graphicx}
\usepackage{color}
\usepackage[pdfstartview=FitH ]{hyperref}
\usepackage[hyperpageref]{backref}%
\newtheorem{theorem}{Theorem}
\newtheorem{lemma}{Lemma}

\newtheorem{proposition}[theorem]{Proposition}

\newenvironment{proof}[1][Proof]{\noindent\textbf{#1.} }{\ \rule{0.5em}{0.5em}}
\begin{document}

\title{The Analytic Strong Multiplicity One Theorem for $\mathrm{GL}_{m}\left(
\mathbb{A}_{K}\right)  $}
\author{Wang Yonghui\\Department of Mathematics\\Capital Normal University\\Beijing 100037 \\P.R. China\\arithwsun@gmail.com or yhwang@mail.cnu.edu.cn}
\date{}
\maketitle

\begin{abstract}
Let $\pi=\otimes\pi_{v}$ and $\pi^{\prime}=\otimes\pi_{v}^{\prime}$ be two
irreducible, automorphic, cuspidal representations of $GL_{m}\left(
\mathbb{A}_{K}\right)  .$ Using the logarithmic zero-free region of
Rankin-Selberg $L$-function, Moreno established the analytic strong
multiplicity one theorem if at least one of them is self-contragredient, i.e.
$\pi$ and $\pi^{\prime}$ will be equal if they have finitely many same local
components $\pi_{v},\pi_{v}^{\prime},$ for which the norm of places are
bounded polynomially by the analytic conductor of these cuspidal
representations. Without the assumption of the self-contragredient for
$\pi,\pi^{\prime},$ Brumley generalized this theorem by a different method,
which can be seen as an invariant of Rankin-Selberg method. In this paper,
influenced by Landau's smooth method of Perron formula, we improved the degree
of Brumley's polynomial bound to be $4m+\varepsilon.$

\end{abstract}

\section{Introduction}

\footnotetext{2000 Mathematics Subject Classification. Primary 11F70, 11S40,
11F67.}\footnotetext{The author is supported by China NSF Grant \#10601034 and
Morning Side center.}

For $K$ an algebraic number field, let $\pi,\pi^{\prime}$ be two cuspidal
automorphic representations of $GL_{m}\left(  \mathbb{A}_{K}\right)  ,$ with
restricted tensor product decompositions $\pi=\otimes\pi_{v}$ and $\pi
^{\prime}=\otimes\pi_{v}^{\prime}$. Proved by Casselman \cite{Casselman1973},
Shalika \cite{Shalika1974}, Piatetski-Shapiro \cite{P-Shapiro1975}, Gelfand
and Kazhdan \cite{GK1975}, Jacquet and Shalika \cite{JS1981}, the strong
multiplicity one theorem states that if $\pi_{v}\cong\pi_{v}^{\prime}$ for all
but finitely many places $v,$ then $\pi\cong\pi^{\prime}.$

By exploiting the logarithmic zero-free region of the Rankin-Selberg
$L$-function $L\left(  s,\pi\times\tilde{\pi}^{\prime}\right)  ,$ Moreno
\cite{Moreno1985} proved the following analytic form of strong multiplicity
one theorem:

\begin{theorem}
[Moreno]\label{Moreno} Let $\mathcal{A}_{n}^{\prime}(Q)$ denote the set of all
cuspidal automorphic representations on $GL_{n}(\mathbb{A}_{K})$ with analytic
conductor less than $Q.$ Suppose $\pi=\otimes\pi_{v}$ and $\pi^{\prime
}=\otimes\pi_{v}^{\prime}$ are in $\mathcal{A}_{n}^{\prime}(Q),$ $n\geq2.$
Then there exist positive constants $C$ and $A_{1}$ such that, if $\pi
_{v}\cong\pi_{v}^{\prime}$ for all finite places $v$ with norm $N\left(
v\right)  \leq Cf\left(  Q\right)  $, then $\pi\cong\pi^{\prime}.$ Here
$f\left(  Q\right)  =Q^{A_{1}}$ for $n=2$ and $f\left(  Q\right)
=e^{A_{1}\left(  \log Q\right)  ^{2}}$ for $n\geq3.$
\end{theorem}

The bound of exponential type in the case of $n\geq3$ is rather poor. The
reason is: until now, we are only able to detect the logarithmic zero-free
region for the Rankin-Selberg $L$-function $L\left(  s,\pi\times\tilde{\pi
}^{\prime}\right)  $ while at least one of $\pi,$ $\pi^{\prime}$ is
self-contragredient. This can be done by applying de la Vall\'{e}e Poussin
method to the corresponding positive $L$-function of the following isobaric
automorphic representations \cite{WB}:%
\[
\Pi=\pi\boxplus\left(  \pi\otimes\alpha^{it}\right)  \boxplus\left(
\pi\otimes\alpha^{-it}\right)  \boxplus\pi^{\prime}\boxplus\left(  \pi
^{\prime}\otimes\alpha^{it}\right)  \boxplus\left(  \pi^{\prime}\otimes
\alpha^{-it}\right)  ,
\]
if $\pi,$ $\tilde{\pi}^{\prime}$ are both self-contragredient (see Sarnak
\cite{Sarnak2003}), and to%
\[
\Pi=\pi\boxplus\left(  \pi^{\prime}\otimes\alpha^{-it}\right)  \boxplus\left(
\tilde{\pi}^{\prime}\otimes\alpha^{it}\right)  ,
\]
if $\pi$ is self-contragredient and $\tilde{\pi}^{\prime}$ is not. Hence, if
at least one of $\pi,$ $\tilde{\pi}^{\prime}$ is self-contragredient, by the
same argument of Moreno (see \cite[p. 183, footnote]{Moreno1985}), $f\left(
Q\right)  =Q^{A_{1}}$ also holds for $n\geq3.$

If both $\pi,$ $\tilde{\pi}^{\prime}$ are not self-contragredient for $n\geq
3$, it is only known that $L\left(  s,\pi\times\tilde{\pi}^{\prime}\right)  $
is nonvanishing for $\operatorname{Re}s\geq1$. And then only the exponential
bound can be obtained under Moreno's method. Recently, Brumley
\cite{Brumley2004} used a different method without applying the derivative of
logarithmic of $L\left(  s,\pi\times\tilde{\pi}^{\prime}\right)  $ as Moreno
did, and obtained that

\begin{theorem}
[Brumley]Suppose $\pi=\otimes\pi_{v}$ and $\pi^{\prime}=\otimes\pi_{v}%
^{\prime}$ are in $\mathcal{A}_{n}^{\prime}(Q),$ $n\geq1.$ Denote by $S$ the
set of all finite places of $K$ at which either $\pi_{v}$ or $\pi_{v}^{\prime
}$ is ramified. There exists positive constant $C=C\left(  n\right)  $ and
$A_{1}=A_{1}\left(  n\right)  $ such that, if $\pi_{v}\cong\pi_{v}^{\prime}$
for all finite places $v\notin S$ with norm $N\left(  v\right)  \leq
CQ^{A_{1}},$ then $\pi\cong\pi^{\prime}.$
\end{theorem}

In this paper, we reorganize Brumley's method and use a appropriate Perron
formula, which can be traced back to Landau's method \cite{Landau1915} on the
$L$-functions with multiple gamma factors. We improve Brumley's result as following:

\begin{theorem}
\label{thm:M-1}Let $\mathcal{A}_{N}(Q)$ denote the set of all cuspidal
automorphic representations on $GL_{m}(\mathbb{A}_{K})$ $\left(  1\leq m\leq
N\right)  $ with analytic conductor less than $Q.$ If $\pi=\otimes\pi_{v}$ and
$\pi^{\prime}=\otimes\pi_{v}^{\prime}$ are in $\mathcal{A}_{N}(Q)$. Then there
exists a constant $C\left(  \varepsilon\right)  $ depending only on
$\varepsilon>0$, $K$ and $N\ $such that, If $\pi_{v}\simeq\pi_{v}^{\prime}$
for all finite places with norm $N\left(  v\right)  <C\left(  \varepsilon
\right)  Q^{4N+\varepsilon}$, then $\pi\cong\pi^{\prime}.$
\end{theorem}

We notice that, applying the Riemann-Roch theorem on the modular curve
$X_{0}\left(  M\right)  ,$ Murty \cite{Murty1997} showed that when $\pi$ and
$\pi^{\prime}$ correspond to holomorphic modular forms of level $M$ and even
weight $k$, then $\pi\cong\pi^{\prime}$ if $\pi_{p}\cong\pi_{p}^{\prime}$ for
all $p<CkM\log\log M$, which means $A_{1}\left(  N\right)  =1+\varepsilon$ in
that case.

\section{Preliminaries on Automorphic $L$-function}

If $\pi$ is an automorphic irreducible cuspidal representation of
$GL_{m}\left(  \mathbb{A}_{K}\right)  ,$ $K$ is the algebraic number field of
degree $[K:\mathbb{Q}]=l=r_{1}+2r_{2},$ $d_{K}$ its discriminant (in the
standard notation$).$ Then (see \cite[p. 47]{Cogdell}) at every finite place
$v$ where $\pi_{v}$ is unramified we have associated a semisimple conjugacy
class, say
\[
A_{\pi,v}=\left(
\begin{array}
[c]{ccc}%
\alpha_{\pi,v}\left(  1\right)  &  & \\
& \ddots & \\
&  & \alpha_{\pi,v}\left(  m\right)
\end{array}
\right)  ,
\]
so that the local $L$-factors are defined by%
\begin{equation}
L\left(  s,\pi_{v}\right)  =\det\left(  I-q_{v}^{-s}A_{\pi,v}\right)
^{-1}=\prod_{i=1}^{m}\left(  1-\alpha_{\pi,v}\left(  i\right)  q_{v}%
^{-s}\right)  ^{-1} \label{eq:Euler-1}%
\end{equation}
where $q_{v}=|\bar{\omega}_{v}|_{v}^{-1}=N\left(  \mathfrak{p}_{v}\right)
=N\left(  v\right)  $ is the module of $K_{v}.$ For the other finite places
$v,$ it is possible to write the local factors at ramified places in the form
of (\ref{eq:Euler-1}) with the convention that some of the $\alpha_{\pi
,v}\left(  i\right)  $'s can be zero. For the infinite places, we write as
\cite[(2.2)]{RW2003}
\[
L\left(  s,\pi_{\infty}\right)  =\pi^{-lms/2}\prod_{i=1}^{lm}\Gamma\left(
\frac{s+b_{\pi}\left(  i\right)  }{2}\right)  ,
\]
Denote by%
\begin{align*}
L\left(  s,\pi_{f}\right)   &  =\prod_{v<\infty}L\left(  s,\pi_{v}\right)
\text{ for }\operatorname{Re}s>>1,\\
L\left(  s,\pi\right)   &  =L\left(  s,\pi_{\infty}\right)  L\left(  s,\pi
_{f}\right)  \text{ for }\operatorname{Re}s>>1,
\end{align*}
in which the absolute convergence can be provided by the following work of
Luo-Rudnick-Sarnak \cite{LRS},%
\begin{equation}
\left\vert \log_{N(v)}\left\vert \alpha_{\pi,v}\left(  i\right)  \right\vert
\right\vert ,\text{ }\operatorname{Re}b_{\pi}\left(  i\right)  \leq
1/2-1/\left(  n^{2}+1\right)  . \label{bound-ramanujan}%
\end{equation}
$L\left(  s,\pi\right)  $ has an analytic continuation and entire everywhere,
we also have the functional equation that%
\[
L\left(  s,\pi\right)  =W_{\pi}q_{\pi}^{\frac{1}{2}-s}L\left(  1-s,\tilde{\pi
}\right)
\]
where $\tilde{\pi}$ is the contragredient of $\pi,$ $W_{\pi}$ is a complex
constant, and $q_{\pi}>0$ is called $``$arithmetic conductor$"$.

For $\pi$ and $\pi^{\prime}$ are automorphic irreducible cuspidal
representations of $GL_{m}\left(  \mathbb{A}_{K}\right)  $ and $GL_{m\prime
}\left(  \mathbb{A}_{K}\right)  $ respectively. The associated Rankin-Selberg
$L$-function is given as an Euler product of degree $mm^{\prime}$%
\[
L\left(  s,\pi\times\tilde{\pi}^{\prime}\right)  =\prod_{v}L_{v}\left(
s,\pi_{v}\times\tilde{\pi}_{v}^{\prime}\right)  \text{ for }\operatorname{Re}%
s>>1.
\]
We can write it as%
\begin{align*}
L\left(  s,\pi\times\tilde{\pi}^{\prime}\right)   &  =L\left(  s,\pi_{\infty
}\times\tilde{\pi}_{\infty}^{\prime}\right)  L\left(  s,\pi_{f}\times
\tilde{\pi}_{f}^{\prime}\right)  ,\\
L\left(  s,\pi_{\infty}\times\tilde{\pi}_{\infty}^{\prime}\right)   &
=\pi^{-mm^{\prime}ls/2}\prod_{i=1}^{mm^{\prime}l}\Gamma\left(  \frac
{s+b_{\pi,\tilde{\pi}^{\prime}}\left(  i\right)  }{2}\right)  ,\\
L\left(  s,\pi_{f}\times\tilde{\pi}_{f}^{\prime}\right)   &  =\prod_{v<\infty
}L_{v}\left(  s,\pi_{v}\times\tilde{\pi}_{v}^{\prime}\right)  \text{ for
}\operatorname{Re}s>>1.
\end{align*}
And%
\begin{equation}
L\left(  s,\pi_{v}\times\tilde{\pi}_{v}^{\prime}\right)  =\prod_{i=1}^{m}%
\prod_{j=1}^{m^{\prime}}\left(  1-\alpha_{\pi,v}\left(  i\right)  \bar{\alpha
}_{\pi^{\prime},v}\left(  j\right)  q_{v}^{-s}\right)  ^{-1}
\label{eq:coef-unramified}%
\end{equation}
are finite local $L$-factors for unramified finite places $v$ (i.e. $\pi_{v}$
and $\pi_{v}^{\prime}$ are both unramified). One can consult \cite[p.
36]{Cogdell} on the discussions for the other ramified cases. Anyway,
$L\left(  s,\pi_{f}\times\tilde{\pi}_{f}^{\prime}\right)  $ defines a
Dirichlet series:%
\begin{equation}
L\left(  s,\pi_{f}\times\tilde{\pi}_{f}^{\prime}\right)  =\sum_{n=1}^{\infty
}a_{\pi\times\tilde{\pi}^{\prime}}\left(  n\right)  n^{-s}. \label{eq:coef-a}%
\end{equation}
The theory of $L\left(  s,\pi\times\tilde{\pi}^{\prime}\right)  ,L\left(
s,\pi_{f}\times\tilde{\pi}_{f}^{\prime}\right)  $ have been developed by
Jacquet, Piatetski-Shapiro, Shalika, Shahidi, Moeglin and Waldspurger:

\begin{enumerate}
\item[RS1] The Euler product for $L\left(  s,\pi_{f}\times\tilde{\pi}%
_{f}^{\prime}\right)  $ converges absolutely for $\operatorname{Re}s>1$
(Jacquet and Shalika \cite{JS1981}). $L(s,\pi\times\tilde{\pi}^{\prime})$ and
$L(s,\pi_{f}\times\tilde{\pi}_{f}^{\prime})$ are non-zero in
$\operatorname{Re}s\geq1.$ (Shahidi \cite{Sha1})

\item[RS2] The complete $L$-function $L\left(  s,\pi\times\tilde{\pi}^{\prime
}\right)  $ has an analytic continuation to the entire complex plane and
satisfies a functional equation%
\[
L\left(  s,\pi\times\tilde{\pi}^{\prime}\right)  =W_{\pi\times\tilde{\pi
}^{\prime}}q_{\pi\times\tilde{\pi}^{\prime}}^{\frac{1}{2}-s}L\left(
1-s,\tilde{\pi}\times\pi^{\prime}\right)
\]
where the $``$root number$"$ $W_{\pi\times\tilde{\pi}^{\prime}}$ is a complex
constant of modulus $1$ and $q_{\pi\times\tilde{\pi}^{\prime}}=d_{K}%
^{mm^{\prime}}N_{\pi\times\tilde{\pi}^{\prime}}>0$ is the "\textit{arithmetic
conductor}". (Shahidi \cite{Sha1}, \cite{Sha2}, \cite{Sha3}, \cite{Sha4}).

\item[RS3] Denote by $\alpha(g)=|\det(g)|$. When $\pi^{\prime}\not \cong
\pi\otimes\alpha^{it}$ for any $t\in\mathbb{R}$, $L(s,\pi\times\tilde{\pi
}^{\prime})$ is holomorphic. When $m=m^{\prime}$ and $\pi^{\prime}\cong
\pi\otimes\alpha^{i\tau_{0}}$ for some $\tau_{0}\in\mathbb{R}$, the only poles
of $L(s,\pi\times\tilde{\pi}^{\prime})$ are simple poles at $s=i\tau_{0}$ and
$1+i\tau_{0}$. (Jacquet and Shalika \cite{JS1981} and Moeglin and Waldspurger
\cite{MoeWa1989}).

\item[RS4] $L\left(  s,\pi\times\tilde{\pi}^{\prime}\right)  $ is meromorphic
of order one away from its poles, and bounded in the vertical strips (see
\cite{GS2001}).
\end{enumerate}

Denote by the $``$\textit{analytic conductor}$"$ for $L\left(  s,\pi\right)  $
and $L\left(  s,\pi\times\pi^{\prime}\right)  $ as%
\begin{align*}
C\left(  \pi\right)   &  =q_{\pi}\prod_{i=1}^{ml}\left(  1+\left\vert b_{\pi
}\left(  i\right)  \right\vert \right) \\
C\left(  \pi,\pi^{\prime}\right)   &  =q_{\pi\times\pi^{\prime}}\prod
_{i=1}^{mm^{\prime}l}\left(  1+\left\vert b_{\pi\times\tilde{\pi}^{\prime}%
}\left(  i\right)  \right\vert \right)  .
\end{align*}
Note that the definition here differs only slightly from that given in Iwaniec
and Sarnak's article \cite{IS2000}, in which we ignore the imaginary
parameter. Ramakrishan and S. Wang named such definition as $``$thickened
conductor$"$ and proved that \cite[(2.8)]{RW2003}%
\begin{equation}
C\left(  \pi\right)  ^{-m^{\prime}}C\left(  \pi^{\prime}\right)  ^{-m}\ll
C\left(  \pi,\pi^{\prime}\right)  \ll C\left(  \pi\right)  ^{m^{\prime}%
}C\left(  \pi^{\prime}\right)  ^{m}. \label{eq:Bound-cond}%
\end{equation}

\section{Proof of Theorem \ref{thm:M-1}}

In the following proof, without loss of generality, we can assume that
$\pi,\pi^{\prime}$ is not twisted-equivalent, i.e. $\pi^{\prime}\ncong
\pi\otimes\alpha^{i\tau_{0}}$ for any $\tau_{0}\neq0.$ Since if $\pi^{\prime
}\cong\pi\otimes\alpha^{i\tau_{0}}$ with $\tau_{0}\neq0,$ then $\pi_{v}%
\ncong\pi_{v}^{\prime}$ for many finite places $v.$ We can write the finite
part of the Rankin-Selberg $L$-function as an Dirichlet series:%
\[
L\left(  s,\pi_{f}\times\tilde{\pi}_{f}^{\prime}\right)  =\prod_{v<\infty
}L_{v}\left(  s,\pi_{v}\times\tilde{\pi}_{v}^{\prime}\right)  =\sum
_{n=1}^{\infty}\frac{a_{\pi\times\tilde{\pi}^{\prime}}\left(  n\right)
}{n^{s}},\text{ }\operatorname{Re}s>1.
\]

Denote by%
\[
S\left(  x;\pi,\tilde{\pi}^{\prime}\right)  :=\sum_{n=1}^{\infty}a_{\pi
\times\tilde{\pi}^{\prime}}\left(  n\right)  w\left(  \frac{n}{x}\right)  ,
\]
where $w(x)$ is a nonnegative real-value function of $C_{c}^{\infty}$ with
compact support in $[0,3]$ and satisfies%
\[
w\left(  x\right)  =\left\{
\begin{array}
[c]{l}%
0,\text{ for }x\leq0\text{ or }x\geq3,\\
e^{-\frac{1}{x}},\text{ for }0<x\leq1,\\
e^{-\frac{1}{3-x}},\text{ for }2<x<3.
\end{array}
\right.
\]
Hence for arbitrary integer $k$, the derivative $w^{(k)}\left(  x\right)  $ is
also exponential decay as $x\rightarrow0,3.$

We first need a Lemma of Brumley \cite{Brumley2004} to estimate the lower
bound of $S\left(  x;\pi,\tilde{\pi}^{\prime}\right)  $ for $\pi^{\prime}%
=\pi.$ Here, we call $a_{\pi\times\tilde{\pi}^{\prime}}\left(  n\right)  $
unramified coefficients if it is only determined by the product of unramified
$L\left(  s,\pi_{v}\times\tilde{\pi}_{v}^{\prime}\right)  ,$ (see
(\ref{eq:coef-unramified})).

\begin{lemma}
[Brumley]If $\pi$ is an automorphic irreducible cuspidal representation of
$GL_{m}\left(  \mathbb{A}_{K}\right)  ,$ $m\geq1.$ $\mathfrak{p}$ denotes any
unramified places (prime ideals),\ then the unramified coefficients for its
Rankin-Selberg $L$-function
\[
a_{\pi\times\tilde{\pi}}\left(  \mathfrak{p}^{m}\right)  \geq1.
\]

\end{lemma}

\begin{proof}
See Lemma 1 of Brumley \cite{Brumley2004}.
\end{proof}

Hence, we can obtain that

\begin{proposition}
If $\pi$ is an automorphic irreducible cuspidal representation of
$GL_{m}\left(  \mathbb{A}_{K}\right)  ,$ then there exist a positive constant
$c=c\left(  K\right)  $ only dependent of the number field $K,$ such that%
\begin{equation}
S\left(  x;\pi,\tilde{\pi}\right)  \geq\frac{c}{2e^{2}}\frac{x^{1/m}}{\log x}
\label{eq:mainterm}%
\end{equation}
for $x\geq\left(  \log Q\right)  ^{3m}.$
\end{proposition}

\begin{proof}
As the coefficients $a_{\pi\times\tilde{\pi}}\left(  n\right)  $ are
nonnegative, the sum $S\left(  x;\pi,\tilde{\pi}\right)  $ can be truncated to
give
\begin{align*}
S\left(  x;\pi,\tilde{\pi}\right)   &  =\sum_{n=1}^{\infty}a_{\pi\times
\tilde{\pi}}\left(  n\right)  w\left(  \frac{n}{x}\right)  \geq\sum_{x/2\leq
n\leq x}a_{\pi\times\tilde{\pi}}\left(  n\right)  e^{-x/n}\\
&  \geq e^{-2}\sum_{x/2\leq n\leq x}a_{\pi\times\tilde{\pi}}\left(  n\right)
\geq e^{-2}\sum_{x/2\leq N\mathfrak{p}^{m}\leq x}a_{\pi\times\tilde{\pi}%
}\left(  \mathfrak{p}^{m}\right) \\
&  \geq e^{-2}\sum_{\substack{\left(  x/2\right)  ^{1/m}\leq N\mathfrak{p}\leq
x^{1/m}\\\mathfrak{p}\text{ is unramified}}}1.
\end{align*}
Therefore, using the prime ideals number theorem, and $\#\left\{
\mathfrak{p}\mid\mathfrak{p}\text{ is ramified}\right\}  \leq\log Q,$ we get%
\begin{align*}
S\left(  x;\pi,\tilde{\pi}\right)   &  \geq e^{-2}\#\left\{  \mathfrak{p}%
\mid\left(  x/2\right)  ^{1/m}\leq N\mathfrak{p}\leq x^{1/m},\mathfrak{p}%
\text{ is unramified}\right\} \\
&  \geq e^{-2}\#\left\{  \mathfrak{p}\mid\left(  x/2\right)  ^{1/m}\leq
N\mathfrak{p}\leq x^{1/m}\right\} \\
&  -e^{-2}\#\left\{  \mathfrak{p}\mid\left(  x/2\right)  ^{1/m}\leq
N\mathfrak{p}\leq x^{1/m},\mathfrak{p}\text{ is ramified}\right\} \\
&  \geq\frac{1}{2e^{2}}\left\{  \mathfrak{p}\mid\left(  x/2\right)  ^{1/m}\leq
N\mathfrak{p}\leq x^{1/m}\right\} \\
&  \geq\frac{c}{2e^{2}}\frac{x^{1/m}}{\log x},
\end{align*}
for $x\geq\left(  \log Q\right)  ^{3m}$ and some\ positive constant
$c=c\left(  K\right)  .$
\end{proof}

Second, let $\mathcal{A}_{N}(Q)$ denote the set of all cuspidal automorphic
representations on $GL_{m}(\mathbb{A})$ $\left(  m\leq N\right)  $ with
analytic conductor less than $Q.$ We assume that $\pi,\pi^{\prime}%
\in\mathcal{A}_{N}\left(  Q\right)  $ and $\pi^{\prime}\ncong\pi\otimes
\alpha^{i\tau_{0}}$ for any $\tau_{0}\in\mathbb{R}.$ Hence, $L\left(
s,\pi\times\tilde{\pi}^{\prime}\right)  $ and $L\left(  s,\pi_{f}\times
\tilde{\pi}_{f}^{\prime}\right)  $ are both holomorphic functions. We will
provide an upper bound for $S\left(  x;\pi,\tilde{\pi}^{\prime}\right)  $ by
the classical analytic argument, which originates from Landau's method
\cite{Landau1915}.

Let%
\begin{align*}
G\left(  s\right)   &  :=\frac{L\left(  1-s,\tilde{\pi}_{\infty}\times
\pi_{\infty}^{\prime}\right)  }{L\left(  s,\pi_{\infty}\times\tilde{\pi
}_{\infty}^{\prime}\right)  }\\
&  =\pi^{-\frac{mm^{\prime}l}{2}+mm^{\prime}ls}\prod_{i=1}^{mm^{\prime}l}%
\frac{\Gamma\left(  \left(  1-s+b_{\tilde{\pi},\pi^{\prime}}\left(  i\right)
\right)  /2\right)  }{\Gamma\left(  (s+b_{\pi,\tilde{\pi}^{\prime}}\left(
i\right)  )/2\right)  }\\
&  =\pi^{-\frac{mm^{\prime}l}{2}+mm^{\prime}ls}\prod_{i=1}^{mm^{\prime}l}%
\frac{\Gamma\left(  \left(  1-s+\bar{b}_{\pi,\tilde{\pi}^{\prime}}\left(
i\right)  )\right)  /2\right)  }{\Gamma\left(  (s+b_{\pi,\tilde{\pi}^{\prime}%
}\left(  i\right)  /2\right)  }.
\end{align*}
For every fixed strip $\sigma_{1}\leq\sigma\leq\sigma_{2}$ and uniformly
growing positive $t>1,$ it is known by Stirling formula that%
\begin{align*}
\Gamma\left(  \sigma+it\right)   &  =c\left(  \sigma\right)  e^{-\frac{\pi}%
{2}t}t^{\sigma-\frac{1}{2}}e^{it\left(  \log t-1\right)  }\left(  1+O\left(
\frac{1}{t}\right)  \right) \\
\Gamma\left(  \sigma-it\right)   &  =\bar{c}\left(  \sigma\right)
e^{-\frac{\pi}{2}t}t^{\sigma-\frac{1}{2}}e^{-it\left(  \log t-1\right)
}\left(  1+O\left(  \frac{1}{t}\right)  \right)  .
\end{align*}
Therefore, let $s=\sigma+it,$ $b_{\pi,\tilde{\pi}^{\prime}}\left(  i\right)
=u\left(  i\right)  +\sqrt{-1}v\left(  i\right)  $ i.e. $u\left(  i\right)
=\operatorname{Re}b_{\pi,\tilde{\pi}^{\prime}}\left(  i\right)  ,$ $v\left(
i\right)  =\operatorname{Im}b_{\pi,\tilde{\pi}^{\prime}}\left(  i\right)  ,$
we obtain
\begin{align*}
G\left(  s\right)   &  \ll_{\sigma,N,K}\prod_{i=1}^{mm^{\prime}l}%
\frac{\left\vert t+v\left(  i\right)  \right\vert ^{\frac{1-\sigma+u\left(
i\right)  }{2}-\frac{1}{2}}}{\left\vert t+v\left(  i\right)  \right\vert
^{\frac{\sigma+u\left(  i\right)  }{2}-\frac{1}{2}}}\\
&  \ll_{\sigma,N,K}\prod_{i=1}^{mm^{\prime}l}\left\vert t+v\left(  i\right)
\right\vert ^{\left(  1/2-\sigma\right)  },
\end{align*}
for $t\notin S,$ $S:=\cup_{i}\left\{  t\mid\left\vert t+v\left(  i\right)
\right\vert \leq1\right\}  .$ It is apparent that $S\sqsubseteq\left[
-Q^{2N}-1,Q^{2N}+1\right]  $ since
\[
\left\vert v\left(  i\right)  \right\vert <\prod_{i=1}^{mm^{\prime}l}\left(
1+\left\vert v\left(  i\right)  \right\vert \right)  \leq C\left(  \pi
,\tilde{\pi}^{\prime}\right)  \ll Q^{m+m^{\prime}}%
\]
by formula (\ref{eq:Bound-cond}). Hence, for $1/2-\sigma>0,$ we have
\begin{align}
G\left(  s\right)   &  \ll_{\sigma,N,K}\left(  1+\left\vert t\right\vert
\right)  ^{mm^{\prime}l\left(  1/2-\sigma\right)  }\prod_{i=1}^{mm^{\prime}%
l}\left(  1+\left\vert v\left(  i\right)  \right\vert \right)  ^{\left(
1/2-\sigma\right)  }\label{eq:Bound-G}\\
&  \ll_{\sigma,N,K}\left(  1+\left\vert t\right\vert \right)  ^{mm^{\prime
}l\left(  1/2-\sigma\right)  }Q^{2N\left(  1/2-\sigma\right)  }.\nonumber
\end{align}
for $t\notin S.$

Denote the Mellin transform of $w\left(  x\right)  $ by%
\[
W\left(  s\right)  =\int_{0}^{+\infty}w\left(  x\right)  x^{s-1}dx
\]
It is easily seen that $W\left(  s\right)  $ is an analytic function since the
integrand $w\left(  x\right)  x^{s-1}$ decays exponentially as $x\rightarrow0$
for any complex $s=\sigma+it$. Furthermore, we notice that for $\sigma<-1,$
\begin{equation}
W\left(  s\right)  =\int_{0}^{+\infty}w\left(  x\right)  x^{s-1}%
dx\ll_{A,\sigma}(|t|+1)^{-A}. \label{eq:bound-W}%
\end{equation}
for any $A>1$ by repeated partial summation, the parameters $A$ and $\sigma$
are independent of each other. By the Mellin inversion, we have%
\[
w\left(  x\right)  =\frac{1}{2\pi i}\int_{2-i\infty}^{2+i\infty}W\left(
s\right)  x^{-s}ds.
\]
Hence we obtain the Perron summation formula as following,%
\begin{align*}
S\left(  x;\pi,\tilde{\pi}^{\prime}\right)   &  =\sum_{n=1}^{\infty}%
a_{\pi\times\tilde{\pi}^{\prime}}\left(  n\right)  w\left(  \frac{n}{x}\right)
\\
&  =\frac{1}{2\pi i}\sum_{n=1}^{\infty}a_{\pi\times\tilde{\pi}^{\prime}%
}\left(  n\right)  \int_{2-i\infty}^{2+i\infty}W\left(  s\right)  \left(
\frac{n}{x}\right)  ^{-s}ds\\
&  =\frac{1}{2\pi i}\int_{2-i\infty}^{2+i\infty}x^{s}W\left(  s\right)
L\left(  s,\pi_{f}\times\tilde{\pi}_{f}^{\prime}\right)  ds,
\end{align*}
where the interchange of summation and integral is provided by the absolute
convergence on the line $\operatorname{Re}s=2.$

The polynomial convexity bound of $L\left(  s,\pi_{f}\times\tilde{\pi}%
_{f}^{\prime}\right)  $ can be easily obtained from RS1, RS2 and
(\ref{eq:Bound-G}) by Phragm\'{e}n-Lindel\"{o}f principle. Noticing the rapid
decay of $W\left(  s\right)  $ (\ref{eq:bound-W}), \ we can apply Cauchy
theorem and RS3, i.e. $L\left(  s,\pi_{f}\times\tilde{\pi}_{f}^{\prime
}\right)  $ has no pole, to obtain that
\[
S\left(  x;\pi,\tilde{\pi}^{\prime}\right)  =\frac{1}{2\pi i}\int_{-H-i\infty
}^{-H+i\infty}x^{s}W\left(  s\right)  L\left(  s,\pi_{f}\times\tilde{\pi}%
_{f}^{\prime}\right)  ds.
\]
for any $H>0$ satisfying that $G\left(  s\right)  $ does not have poles at the
line $\operatorname{Re}s=-H.$ Using the functional equation, we have%
\begin{align}
S\left(  x;\pi,\tilde{\pi}^{\prime}\right)   &  =\frac{1}{2\pi i}\int_{\left(
-H\right)  }x^{s}W\left(  s\right)  W_{\pi\times\tilde{\pi}^{\prime}}%
q_{\pi\times\tilde{\pi}^{\prime}}^{1/2-s}G\left(  s\right)  L\left(
1-s,\tilde{\pi}_{f}\times\pi_{f}^{\prime}\right)  ds\nonumber\\
&  =\frac{1}{2\pi i}\int_{\left(  -H\right)  }x^{s}W\left(  s\right)
W_{\pi\times\tilde{\pi}^{\prime}}q_{\pi\times\tilde{\pi}^{\prime}}%
^{1/2-s}G\left(  s\right)  \left(  \sum_{n=1}^{\infty}\frac{a_{\tilde{\pi
}\times\pi^{\prime}}\left(  n\right)  }{n^{1-s}}\right)  ds\nonumber\\
&  =\sum_{n=1}^{\infty}\frac{a_{\tilde{\pi}\times\pi^{\prime}}\left(
n\right)  }{n^{1+H}}\int_{\left(  -H\right)  }x^{s}W\left(  s\right)
W_{\pi\times\tilde{\pi}^{\prime}}q_{\pi\times\tilde{\pi}^{\prime}}%
^{1/2-s}n^{it}G\left(  s\right)  ds, \label{eq:S-1}%
\end{align}
where we can interchange the integral and the summation by the absolute
convergence of the Dirichlet series and the rapid decay of $W\left(  s\right)
$ (\ref{eq:bound-W}). Hence, by the absolute convergence of the Dirichlet
series $L\left(  s,\tilde{\pi}_{f}\times\pi_{f}^{\prime}\right)  $ for
$\operatorname{Re}s>1$ and the upper bound (\ref{eq:Bound-G}) of $G\left(
s\right)  ,$ $S\left(  x;\pi,\tilde{\pi}^{\prime}\right)  $ becomes%
\begin{align*}
&  \sum_{n=1}^{\infty}\frac{a_{\tilde{\pi}\times\pi^{\prime}}\left(  n\right)
}{n^{1+H}}\int_{\left(  -H\right)  }x^{s}W\left(  s\right)  W_{\pi\times
\tilde{\pi}^{\prime}}q_{\pi\times\tilde{\pi}^{\prime}}^{1/2-s}n^{it}G\left(
s\right)  ds\\
&  \ll_{H,N,K}\int_{\left(  -H\right)  }\left\vert x^{s}W\left(  s\right)
W_{\pi\times\tilde{\pi}^{\prime}}q_{\pi\times\tilde{\pi}^{\prime}}%
^{1/2-s}n^{it}G\left(  s\right)  \right\vert ds=\int_{t\in S}+\int_{t\notin
S}\\
&  \ll_{H,N,K}Q^{4N\left(  1/2+H\right)  }x^{-H}.
\end{align*}
Hence, we have

\begin{proposition}
If $\pi,\pi^{\prime}\in\mathcal{A}_{N}\left(  Q\right)  $ and $\pi,\pi
^{\prime}\in\mathcal{A}_{N}\left(  Q\right)  $ and $\pi^{\prime}\ncong
\pi\otimes\alpha^{i\tau_{0}}$ for any $\tau_{0}\in\mathbb{R},$ then
\end{proposition}

\begin{equation}
S\left(  x;\pi,\tilde{\pi}^{\prime}\right)  =O_{H,N,K}\left(  Q^{4N\left(
1/2+H\right)  }x^{-H}\right)  , \label{eq:errorterm}%
\end{equation}

Now, we can conclude

\begin{proof}
[Proof of Theorem \ref{thm:M-1}]If $\pi,\pi^{\prime}$ is not
twisted-equivalent, i.e. $\pi^{\prime}\ncong\pi\otimes\alpha^{i\tau_{0}}$ for
any $\tau_{0}\neq0,$ then $\pi_{v}\ncong\pi_{v}^{\prime}$ for many finite
places $v.$ If $\pi,\pi^{\prime}\in\mathcal{A}_{N}\left(  Q\right)  ,$ not
twisted-equivalent and also $\pi\ncong\pi^{\prime}$, combined with
(\ref{eq:mainterm}) and (\ref{eq:errorterm}), we deduce that%
\begin{equation}
S\left(  x;\pi,\tilde{\pi}\right)  -S\left(  x;\pi,\tilde{\pi}^{\prime
}\right)  \geq\frac{c}{2e^{2}}\frac{x^{1/N}}{\log x}+O_{H,N,K}\left(
Q^{4N\left(  1/2+H\right)  }x^{-H}\right)  . \label{eq:Bound-sum}%
\end{equation}
There exists a constant $C=C\left(  H,N,K\right)  $, such that if we take%
\[
x=CQ^{4N+\frac{2N}{H}},
\]
then the main term is greater than the error term in the formula
(\ref{eq:Bound-sum}). Assume that $\pi_{v}\simeq\pi_{v}^{\prime}$ for all
finite places with norm $N\left(  v\right)  <3x,$ \ it is apparent that%
\[
S\left(  x;\pi,\tilde{\pi}^{\prime}\right)  -S\left(  x;\pi,\tilde{\pi
}\right)  =0.
\]
There will be a contradiction by formula (\ref{eq:Bound-sum}) and the choose
of $x$. Therefore, Theorem \ref{thm:M-1} follows by taking $H$ sufficiently large.
\end{proof}

\textbf{Ackowledgment. }The author would like to thank Brumley for presenting
me his final version of preprint. I am also happy to show my gratitude for the
helpful conversation with Jianya Liu, Yangbo Ye, Chaohua Jia and Claus Bauer.

\end{document}